\documentclass[12pt, reqno]{article}

\usepackage{graphics}
\usepackage{amsmath}
\usepackage{amsfonts}
\usepackage{latexsym}
\usepackage{amssymb}

\newtheorem{theorem}{Theorem}

\newtheorem{lemma}[theorem]{Lemma}

\numberwithin{equation}{section}

\newcommand{\proofend}{$\Box$\bigskip}

\newcommand{\Z}{{\mathbf Z}}

\newcommand{\R}{{\mathbf R}}

\newcommand{\kk}{{\mathbf k}}

\newcommand{\Id}{\rm Id}

\title{Periodic trajectories in 3-dimensional convex billiards}
\author{Michael Farber\footnote{Partially supported by the US -
Israel Binational
Science Foundation; part of
this work was
done while M.F. was visiting Max-Planck Institute for Mathematics in
Bonn} \, and Serge Tabachnikov}        

\date{April 5, 2001}          

\begin{document}
\maketitle

\thanks{

\begin{abstract}
We give a lower bound on the number of periodic billiard trajectories
inside a generic smooth
strictly  convex closed surface in 3-space: for odd $n$, there are at
least $2(n-1)$
such trajectories. Convex plane billiards were studied by G.
Birkhoff, and the case of
higher dimensional billiards is considered in our previous papers.
We apply a topological approach based on the calculation of cohomology
of certain configuration spaces of points on 2-sphere.

{\it MSC}: Primary 3Dxx;  Secondary 58Exx


\end{abstract}

\section{\bf Introduction}

Given a smooth strictly  convex closed hypersurface $X^{m}\subset
\R^{m+1}$, consider
the billiard system inside the convex body $T$ (the billiard table),
bounded by $X$.
One views the billiard ball as a point which moves freely inside $T$
with the elastic reflections
in the boundary $X= \partial T$ subject to the familiar law of
geometric optics: the incoming
and the outgoing trajectories lie in the same 2-plane with the normal
at the impact point, and
the angle of reflection equals the angle of incidence. Billiard
trajectories are extrema of the
length functional.

The problem of periodic billiard trajectories is related to the
closed geodesics problem
on Riemannian manifolds.
Unlike the thoroughly studied latter problem (see \cite{Kl}), the
problem of estimating the number
    of periodic trajectories in convex billiards did not receive much
attention. The first results are due to G.
Birkhoff \cite{Bi} who studied periodic billiard trajectories in
convex plane billiards (when $m=1$), known also as
"Birkhoff billiards". The next step was made by I. Babenko \cite{Ba}
about 10 years ago who considered convex
billiards in 3-space. This case is technically the hardest, and the
paper \cite{Ba} contained an error in
cohomological computations (to quote R. Bott \cite{Bo} on the closed
geodesics problem: ``For some reason, the
large literature on this subject teems with mistakes"). The case of
convex billiards in $m$-dimensional spaces
with $m \geq 3$ was studied in our recent paper \cite{FT}. A relevant
result is as follows.

\begin{theorem}[\cite{FT}] \label{thm1} Let $n \geq 3$ be an odd
number and  $X^{m}\subset \R^{m+1},
m \geq 3$, be a generic smooth strictly  convex closed hypersurface.
Then the number of distinct $n$-periodic
billiard trajectories inside $X$ is not less than $(n-1)m$.
\end{theorem}

Note that the dihedral group $D_n$ acts on $n$-periodic billiard
trajectories in a natural way, and the above
theorem, as  well as the results below, estimates the
number of distinct $D_n$-orbits of  $n$-periodic
trajectories.

The techniques of \cite{FT} do not apply to the case
    $m=2$.

A weaker estimate (which however holds for {\it all} strictly convex billiard
tables without any genericity assumptions) is obtained in \cite{F2}:

\begin{theorem}[\cite{F2}]  Let $n \geq 3$ be an odd prime and $X$ be
a smooth strictly  convex closed
surface in 3-space. Then the number of distinct $D_n$-orbits of
$n$-periodic billiard
trajectories inside $X$ is not less than
$(n+1)/2$.
\end{theorem}

The goal of the present paper is to show that the result of Theorem
\ref{thm1} holds for $m=2$ as well
(it does, for $m=1$, as follows from Birkhoff's work \cite{Bi}).
Our main result is as follows.

\begin{theorem}\label{main}
Let $n \geq 3$ be an odd number and let $X$ be a generic smooth strictly
convex closed surface in 3-space.
Then the number of distinct $D_n$-orbits of $n$-periodic billiard
trajectories inside
$X$ is not less than $2(n-1)$.
\end{theorem}

The meaning of the genericity assumption is explained in the next section.

\section{\bf Cyclic configuration spaces}

    Let $X$ be a topological space and $n$ be a positive integer.
Denote by $G(X,n)$ the subspace of the
Cartesian  power  $X^{\times n} = X\times X\times \dots \times X$
consisting of all configurations $(x_1, x_2,
\dots, x_n)$ such that $x_i\ne x_{i+1}$ for $i=1, 2, \dots, n-1$ and
$x_n \ne x_1$. The space
$G(X,n)$ is called {\it the cyclic configuration space}, see \cite{FT}.

The dihedral group $D_n$ acts naturally on the cyclic configuration
space $G(X,n)$.
One may think of $D_n$ as the group of symmetries of a regular
$n$-gon. Identifying $D_n$ with the subgroup
of permutations of the set of vertices $1, 2, \dots, n$, one
describes the action of $D_n$ on $G(X,n)$ as
follows:
$$\tau\cdot (x_1, \dots, x_n) = (x_{\tau(1)}, \dots, x_{\tau(n)})\ \
{\rm for} \ \ \tau \in D_n, \ (x_1,
\dots, x_n)\in G(X,n).$$
The action of $D_n$ is free if and only if $n$ is  prime.

The strategy of the proof of Theorem \ref{main} (as well as two other
theorems stated in the Introduction)
is as follows.  As was already mentioned, $n$-periodic billiard
trajectories inside the billiard table
$T$ are critical points of the perimeter length function $L_X$ on the
space of $n$-gons, inscribed in $X = \partial T$. The space of
inscribed $n$-gons
is the cyclic
configuration space $G(X,n)$. The genericity
assumption in Theorem \ref{main} means that  $L_X$ is a Morse
function of $G(X,n)$. This function is invariant
under the action of the dihedral group $D_n$. Although the space
$G(X,n)$ is not compact, one can show that the
number of critical $D_n$-orbits of $L_X$ is bounded below by the rank
of the $D_n$-equivariant cohomology space
of the cyclic configuration space -- see \cite{FT}, \S 4.

More precisely, we will use Proposition 4.5 from \cite{FT}.
Since
$X$ is topologically the 2-sphere, our goal is to find the dimension of
the graded vector space $H^\ast_{D_n} (G(S^2,n);\Z_2)$. The answer
is given in the next theorem which, therefore, implies Theorem \ref{main}.

\begin{theorem} \label{Betti} For odd $n \geq 3$, the Poincar\'e
polynomial of $H^\ast_{D_n}(G(S^2,n);\Z_2)$
equals
\[(1+t^2)\cdot(1+t+t^2+\dots+t^{n-2}),\]
and the sum of the respective Betti numbers is $2(n-1)$.
\end{theorem}
 
This result complements the computation of the rings
$H^\ast_{D_n}(G(S^m,n);\Z_2)$ for $m \geq 3$  and
odd $n$ in \cite{FT} and the rings  $H^\ast(G(S^m,n);\kk)$ for
various values of $m$ and various fields
$\kk$ in \cite{F2}. However,  the
results in \cite{F2} do not cover
the case $m=2$ and $\kk=\Z_2$. The calculations of \cite{FT},
\cite{F1}, \cite{F2}
are based on a
spectral sequence
computing  the cohomology algebra of the cyclic configuration space
$G(X,n)$ for an arbitrary manifold $X$, developed in \cite{FT}; this
spectral sequence is
similar to the Totaro spectral sequence
\cite{To} for the usual configuration spaces.

\section{\bf Cohomological computations}

In this section we prove Theorem \ref{Betti}. We will use some
results from our previous papers \cite{FT,F1,F2}. Recall that $n$ is
assumed to be odd.

Let $T: G(S^2,n)\to G(S^2,n)$ be the reflection with respect to the first point
\[T(x_1, x_2, \dots, x_n) = (x_1, x_n, x_{n-1}, \dots, x_2), \quad T^2=\Id.\]
Clearly $T$ has no fixed configurations. 
Denote by $G'=G(S^2,n)/T$ the quotient manifold.

Let $\Z_2\subset D_n$ denote the cyclic subgroup of the dihedral group
$D_n$ generated by $T$.
We want to prove that
$H^\ast_{D_n}(G(S^2,n);\Z_2)$ is isomorphic to
$H^\ast_{\Z_2}(G(S^2,n);\Z_2) \simeq H^\ast(G';\Z_2)$,
and we will compute below the latter space.
Consider the
following commutative diagram
\[
\begin{array}{ccc}
E\Z_2\times _{\Z_2} G(S^2,n) & \to & ED_n\times _{D_n} G(S^2,n)\\ \\
\makebox[3cm][r]{$\downarrow \, G(S^2,n)$} & & \makebox[3cm][r]{$\downarrow \, G(S^2,n)$} \\ \\
B\Z_2 & \to & BD_n
\end{array}
\]
and the Serre spectral sequences with $\Z_2$ coefficients of the
the fibrations represented by the vertical arrows (where $EG \to BG$ denotes
the universal bundle of a group $G$).
We obtain a homomorphism of spectral sequences which
induces a homomorphism
\begin{eqnarray}
H^\ast_{D_n}(G(S^2,n);\Z_2) \to H^\ast_{\Z_2}(G(S^2,n);\Z_2). \label{iso}
\end{eqnarray}

\begin{lemma}\label{lemma} The homomorphism (\ref{iso}) is an isomorphism.
\end{lemma}
\noindent
{\bf Proof}.  On the level of $E_2$-terms, one has a homomorphism
\begin{eqnarray}
E_2^{p,q}=H^p(D_n;H^q(G(S^2,n);\Z_2)) \to
  {E'}_2^{p,q}=H^p(\Z_2;H^q(G(S^2,n);\Z_2)).\label{homo}
\end{eqnarray}
We claim that (\ref{homo}) is an isomorphism.
  The claim would follow from Proposition 10.4 of Chapter 3, \cite{Br}
if one shows that the dihedral group $D_n$ acts trivially on the vector space
\begin{eqnarray*}
H^p(\Z_2;H^q(G(S^2,n);\Z_2)).\label{vectorspace}
\end{eqnarray*}
According to Theorem 4 and Remark 3.3 from \cite{FT}, the cyclic group
$Z_n\subset D_n$  acts trivially on $H^q(G(S^2,n);\Z_2)$ and
the dihedral group $D_n$ acts trivially on $H^q(G(S^2,n);\Z_2)$ for $q$
equal $0, 1, n-1,n$. Moreover, for
$1<q< n-1$ the cohomology space $H^q(G(S^2,n);\Z_2)$ is two-dimensional
and $D_n$ acts trivially on a one-dimensional subspace and on its 
quotient space.
It follows, that for any value of $q$ between $1$ and $n-1$ there exist
two possibilities: either $D_n$ acts trivially on $M=H^q(G(S^2,n);\Z_2)$, or
$M$, viewed as a $\Z_2[\Z_2]$-module, is the unique nontrivial extension
$0\to \Z_2\to M\to \Z_2\to 0$. One checks directly, using the definition of
group cohomology, that in the second case $H^p(\Z_2;M)=0$ for all $p>0$.
In the second case $H^0(\Z_2;M)$ is isomorphic to the
subspace of $\Z_2$-invariants in $M$, which clearly coincides with the subspace
of $D_n$-invariants in $M$.
Hence (\ref{homo}) is an isomorphism for all $p$ and $q$.

The comparison theorem for spectral sequences implies an  isomorphism
\[H^\ast_{D_n}(G(S^2,n);\Z_2) \to H^\ast_{\Z_2}(G(S^2,n);\Z_2),\]
as claimed. \proofend

Now we will study the space $H^\ast(G';\Z_2)$.
    Fix a  point $A\in S^2$. Along with the cyclic configuration space
$G(S^2,n)$, consider its
subspace $G_A$ consisting of the configurations $(x_1, x_2, \dots,
x_n)\in G(S^2,n)$ with $x_1=A$.
Clearly, $G_A$ is $T$-invariant; set $G'_A = G_A/T$.

The standard action of $SO(3)$ on $S^2$ by rotations induces an
action of $SO(3)$
on the cyclic configuration
space $G(S^2,n)$
\begin{eqnarray*}
SO(3)\times G(S^2,n)\to G(S^2,n).
\end{eqnarray*}
    The restriction of this map to $G_A$ gives a continuous map
\begin{eqnarray}
p: SO(3)\times G_A \to G(S^2,n).\label{fibration}
\end{eqnarray}
and we claim that it is a fibration.
To prove this we will show that ({\ref{fibration}) is induced from 
the standard fibration
\begin{eqnarray}
q: SO(3) \to S^2,\quad \mbox{where}\quad q(R)=R(A),\quad R\in SO(3), 
\label{hopf}
\end{eqnarray}
by the map
$f: G(S^2,n) \to S^2$ given by $f(x_1, \dots, x_n) =x_1$.  The
total space $E$ of the induced fibration is the space of all pairs
$(R, (x_1, \dots, x_n))$ with $R(A)=x_1$, where $R\in SO(3)$ and
$(x_1, \dots, x_n)\in G(S^2,n)$. The map $(R, (x_1, \dots, 
x_n))\mapsto (R,(A, R^{-1}(x_2),
\dots, R^{-1}(x_n))$
is a homeomorphism $E\to SO(3)\times G_A$. This proves that 
(\ref{fibration}) is the induced
fibration.

Let $c=(A, x_2, \dots, x_n)\in G_A$ be a fixed
configuration. If a rotation $R\in SO(3)$
    preserves $c$ then $R$ is identical (if $n$ is even the
configuration $(A,-A, A, -A, \dots)$ is invariant
under rotations about the axis $(A, -A)$). The fiber $p^{-1}(c)$ is
identified with the set of  pairs
$(R_\phi, R_{-\phi}c)$, where
$R_\phi\in SO(3)$ denotes the rotation through angle $\phi$ about $A$.
Hence the fiber of (\ref{fibration}) is the circle $S^1$. Factorizing
by $T$, one obtains the fibration
\begin{eqnarray}
p': SO(3)\times G'_A \to G'\label{fibration'}
\end{eqnarray}
with fiber $S^1$.  As above, fibration (\ref{fibration'}) is induced from (\ref{hopf}) by the map
$G'\to S^2$ sending orbit of any configuration $(x_1, x_2, \dots, x_n)\in G(S^2,n)$ to $x_1\in S^2$.

Consider the Serre spectral sequence of (\ref{fibration'}) with 
$\Z_2$ coefficients. The cohomology
of the total space
$$H^\ast(SO(3)\times G'_A;\Z_2) \simeq H^\ast(SO(3);\Z_2)\otimes
H^\ast(G_A';\Z_2)$$
is as follows.

It is well known that $H^\ast(SO(3);\Z_2)=\Z_2[v]/(v^4)$ is the truncated polynomial algebra
with a single
1-dimensional generator $v$ satisfying the relation $v^4=0$. The
cohomology algebra $H^\ast(G_A';\Z)$
with integral coefficients was described in Theorem 12 of \cite{F2}. It has
generators
$$\delta_i\in H^{4i}(G_A';\Z)\simeq \Z,\quad\mbox{where}\quad i =1, 2, \dots,$$
and two other generators
$$a\in H^1(G'_A;\Z)\simeq \Z,\quad\mbox{and}\quad
b\in H^3(G'_A;\Z)\simeq \Z_2,$$
which satisfy the following defining relations:
\[
\begin{array}{ll}
\delta_i\delta_j =
\left(
\begin{array}{c}
2i+2j\\
2i
\end{array} \right)
\cdot \delta_{i+j},\quad &
\delta_{[(n+1)/4]}=0,\\
b^2=0 &
2b=0, \\
a^2=0, &ab=0, \\
    \delta_{k}b =0 \quad \mbox{(only for $n=4k+3$).}&
\end{array}
\]
An additive basis of $H^\ast(G'_A;\Z_2)$ is given by  mod $2$ reductions of the
classes $\delta_i$ of degree $4i$ and classes $\delta_ia$ of degree $4i+1$;
moreover, according to the universal coefficient formula, each integral class
   $\delta_ib$ of order 2 produces two classes in $\Z_2$-cohomology,
one of degree $4i+2$ and one of degree $4i+3$.

We claim that the Poincar\`e polynomial of $H^\ast(G'_A;\Z_2)$ equals
\begin{eqnarray}
1+t+t^2+\dots +t^{n-2}.\label{poly}
\end{eqnarray}
To prove this, consider two cases: $n=4k+3$ and $n=4k+1$. In the
first case the mod $2$ reductions of the
generators $\delta_i$ and $\delta_ia$ produce the Poincar\'e polynomial
$$(1+t)(1+t^4+t^8+\dots +t^{4k}),$$
   and the classes $\delta_ib$ contribute
\begin{eqnarray}
(t^2+t^3)(1+t^4+t^8+\dots+t^{4k-4})\label{sum}
\end{eqnarray}
to the Poincar\'e polynomial.
Summing these two polynomials gives (\ref{poly}).
In the second case, when $n=4k+1$, the contribution of the generators
$\delta_i$ and $\delta_ia$ is
$$(1+t)(1+t^4+\dots +t^{4k-4});$$
   the contributions of the classes
$\delta_ib$ is (\ref{sum}), and the total is again (\ref{poly}).

Next, consider the cohomological spectral sequence of the fibration
(\ref{fibration'}). The term $E_2$ consists
of two rows, and the only possibly non-trivial differential is $d_2$.
We claim that this differential
vanishes. Assuming this claim,
the Poincar\`e polynomial of the base equals the Poincar\`e
polynomial of the total space
divided by the Poincar\`e polynomial of the fiber. We obtain  the
Poincar\`e polynomial of
$H^\ast(G';\Z_2)$:
\begin{eqnarray*}
\frac{(1+t+t^2+t^3)\cdot (1+t+t^2+\dots+t^{n-2})}{1+t} \, =\, \\
=\, (1+t^2)\cdot(1+t+t^2+\dots+t^{n-2})
\end{eqnarray*}
as stated in Theorem \ref{Betti}.

It remains to prove that $d_2 =0$. Let $s \in E_2^{0,1} = H^1 (S^1;\Z_2)$ be the generator.
It suffices to show that $d_2 (s) = 0$. The class $s$ is
transgressive, and $d_2 (s)$ is the mod 2 reduction of the Euler class
of fibration (\ref{fibration'}).  Fibration (\ref{fibration'}) is induced from fibration (\ref{hopf}), cf. above.
The latter is the unit tangent bundle
of
the 2-sphere. The Euler number of $S^2$ equals 2, and the Euler class of (\ref{hopf}) vanishes
in $H^2 (S^2;\Z_2)$.
By functoriality of characteristic classes, the mod 2 reduction of the Euler class of
(\ref{fibration'}) is trivial. 

This completes
the proof of Theorem \ref{Betti}.

\bibliographystyle{amsalpha}

\vskip 1cm

Michael Farber,

Department of Mathematics,

Tel Aviv University, 

Ramat Aviv 69978, Israel

farber@math.tau.ac.il

\vskip 1 cm
Serge Tabachnikov

Department of Mathematics,

Penn State University,

University Park, PA 16802, USA

tabachni@math.psu.edu

\end{document}